\documentclass{elsart}

\begin{document}

\begin{frontmatter}
\title{On the Oscillations of Multiplicative Functions Taking Values $\pm 1$}
\author{Ernest S. Croot III}
\ead{ecroot@math.berkeley.edu}
\address{U. C. Berkeley, 1067 Evans Hall, Berkeley, CA  94720}
\begin{abstract}  
For multiplicative functions $f(n)$, which take on the values
$\pm 1$, we show that under certain conditions on $f(n)$, 
for all $x$ sufficiently large, there are at least 
$x\exp(-7(\log\log x)\sqrt{\log x})$ values of $n \leq x$ for which $f(n(n+1)) = -1$.
\end{abstract}
\end{frontmatter}

\section{Introduction}

Given a completely multiplicative function $f(n)$, which takes on the 
values $\pm 1$, and which has the property that
$$
\sum_{n \leq x} f(n) = o(x),
$$
one might wonder how often $f(n)$ changes sign.  That is, for how
many integers $n \leq x$ is $f(n)=-f(n+1)$?

It is a rather simple matter to prove that the number of $n \leq x$ with $f(n)=f(n+1)$ is 
$\gg x$:  By the pigeonhole principle, one can see that there exist
$\delta, \delta' = 0,1,$ or $2$, with $\delta < \delta'$, such that there are 
$\gg x$ integers $j \leq x/2-1$ with $f(2j+\delta)=f(2j+\delta')$.
Then, since $(2j+\delta')/(\delta,\delta') = (2j+\delta)/(\delta,\delta') + 1$,
and since $f((2j+\delta')/(\delta,\delta')) = f((2j+\delta)/(\delta,\delta'))$,
each such $j$ gives rise to an $n \leq x$ with $f(n)=f(n+1)$.  

Despite how similar these two problems seem (counting $n$'s satisfying 
$f(n)=f(n+1)$, and counting $n$'s satisfying $f(n)=-f(n+1)$), it is an incomparably
more difficult problem to show that $f(n)=-f(n+1)$ for $\gg x$
integers $n \leq x$.  Even so, there has been some progress on this question
in the past decade or so.  For instance, Harman, Pintz and Wolke (see \cite{harman2}) 
proved that for $f(n) = \lambda(n) = (-1)^{\Omega(n)}$ (the Liouville function), 
\begin{equation}\label{liouville_bound}
\#\{ n \leq x\ :\ f(n)=-f(n+1)\} > {x \over \log^{7+o(1)}x}.
\end{equation}
In \cite{hildebrand2}, A. Hildebrand proved the following bound for general 
completely multiplicative functions $f(n)=\pm 1$, 
which gives a much better answer for infinitely many $x$:
$$
\limsup_{x \to \infty} {(\log\log x)^4\#\{n \leq x\ :\ f(n)=-f(n+1)\} \over x} > 0.
$$
Perhaps the methods of Hildebrand can be used to replace the ``limsup'' with
a ``lim'', and thereby give a stronger result than our Theorem 
\ref{main_theorem} below.

In this paper, we prove the following result, which generalizes the result of 
Harman, Pintz, and Wolke (although the lower bound we give is not as sharp
as the one they derive for $f(n)=\lambda(n)$):

\begin{thm}\label{main_theorem} 
Suppose that $f(n)$ is a completely multiplicative function which takes
on the values $\pm 1$, and suppose that
\begin{equation}\label{main_assumption}
\lim_{x \to \infty} {1 \over x} \sum_{n \leq x} f(n) = 0.
\end{equation}
Then, for $x$ sufficiently large,
\begin{equation}\label{main_consequence}
\#\{ n \leq x \ :\ f(n) = -f(n+1)\}\ >\ {x \over L(x)^{7}},
\end{equation}
where $L(x)=\exp((\log\log x)\sqrt{\log x})$.
\end{thm}
\bigskip

\noindent {\bf Remark:}  We could perhaps improve the 
$L(x)^{7}$ to $\exp(C\sqrt{\log x})$ for some constant $C$; however, improving the
$\sqrt{\log x}$ to $(\log x)^\epsilon$ for arbitrary $\epsilon > 0$ seems to
require some new ideas.  The bottleneck to obtaining such results, 
using the method in this paper, is the fact that 
Lemma \ref{prop_2_smooth_lemma} only holds with $\beta_2 = \exp((\log x)^{1/2+o(1)})$.  
If one could prove that this Lemma holds for
$\beta_1, \beta_2 < \log^Bx$, for some $B > 0$, then one could prove an estimate
of the same quality as (\ref{liouville_bound}), but with $7$ replaced with 
some constant $C > 0$.

The method of proof of the Theorem is apparently new, and proceeds by showing that 
\begin{equation}\label{oscillate_discuss}
\#\{n \leq x\ :\ f(n) = -f(\lfloor pn/q\rfloor ) \} > {x \over L(x)^6},
\end{equation}
where $p,q \leq \exp(10\sqrt{\log x})$ are some pair of primes with 
$f(p)=f(q)$.  For each such $n$ counted, we will have that 
$f(pn) = -f(q\lfloor pn/q\rfloor)$; and so, since $pn$ and $q \lfloor pn/q \rfloor$ 
lie in $(pn-q,pn]$, we deduce that this interval contains an integer $m$ with
$f(m)=-f(m+1)$. 

We show that if (\ref{oscillate_discuss}) fails to hold for all such primes $p$
and $q$ above, then 
\begin{equation}\label{alpha_bound}
\#\{n \leq x\ :\ f(n) = -f(\lfloor \alpha n \rfloor )\} < {x \over L(x)^3},
\end{equation}
for all real numbers $\alpha$ in a certain short (but not too short) interval
near $1$.  Showing that (\ref{alpha_bound}) cannot hold, for all such $\alpha$
considered, is a relatively simple task, and can be proved by integrating
the left hand side of (\ref{alpha_bound}) over all such $\alpha$, and then showing that
this integral cannot be too small.  
\bigskip

\section{Proof of Theorem \ref{main_theorem}.}\label{main_theorem_section}
\bigskip

For a given $x > 1$, let $I = I(x)$ denote the interval
$[x/(2\beta_1), x/\beta_1]$, where 
$$
\beta_1 = \beta_1(x) = \exp(10\sqrt{\log x}). 
$$  
For a given positive real number $\alpha$, let
$$
\Sigma(I,\alpha) = \#\{ n \in I\ :\ f(n) = -f(\lfloor \alpha n \rfloor ) \}.
$$

Let $R(x)$ be the set of all rational numbers
$$
\alpha = {p_1 p_2 \cdots p_k \over q_1q_2\cdots q_k} \in (1/2,2),
$$
where $k < \beta_2 = \beta_2(x) = L(x)^3$, and
where $p_1,...,p_k,q_1,...,q_k$ are primes $\leq \beta_1$, with
$p_i/q_i \in (1/2,2)$ and $f(p_i)=f(q_i)$.  We note that 
$R(x)$ contains the number $1$.

The proof of the Main Theorem will follow from the following two Propositions:
\bigskip

\begin{prop}\label{prop_1}  
If $x$ is sufficiently large, and if there exists $\alpha' \in R(x)$ such that
\begin{equation}\label{translation_inequality2}
\Sigma(I, \alpha')\ \geq\ {|I| \over \beta_2}, 
\end{equation}
then (\ref{main_consequence}) holds.
\end{prop}

\begin{prop}\label{prop_2}  For all $x$ sufficiently large,
there exists $\gamma \in R(x)$ such that 
$$
\Sigma(I,\gamma) > {|I| \over 3}.
$$
\end{prop} 

Let $\gamma$ satisfy the conclusion to this last proposition (we assume $x$ is
sufficiently large).  Thus, (\ref{translation_inequality2}) holds for 
$\alpha' = \gamma$.  By Proposition \ref{prop_1}, 
(\ref{main_consequence}) holds, and Theorem \ref{main_theorem} is
proved.   \qed
\bigskip

\section{Proof of Proposition \ref{prop_1}.}

We will prove the contrapositive of this Proposition.  So, suppose that 
(\ref{main_consequence}) fails to hold.  We will show that the hypothesis of
the Proposition is false.

In this proof and in later results we will need the following 
Lemma and its corollaries:

\begin{lem} \label{shift_lemma}
If (\ref{main_consequence}) fails to hold (for a particular value of 
$x$), and if $\delta(n)$ is any integer-valued function of $n$ 
satisfying $|\delta(n)| \leq B < x$, then
\begin{equation}\label{shift_equation}
\sum_{B < n \leq x-B} |f(n + \delta(n))\ -\ f(n)|\ \leq\ {4Bx \over L(x)^{7}}.
\end{equation}
\end{lem}

\begin{pf}
For each $n \in (B,x-B]$ 
where $f(n+\delta(n)) = -f(n)$, the interval $[n,n+\delta(n)]$ (or 
$[n+\delta(n),n]$ if $\delta(n) < 0$) contains consecutive integers $m,m+1$
such that $f(m)=-f(m+1)$.  Thus, the left hand side of (\ref{shift_equation})
is bounded from above by 
$$
2 \sum_{m \leq x-1 \atop f(m)=-f(m+1)} \#\{ n \leq x\ :\ |n-m| \leq |\delta(n)|\}
\ \leq\ 2 \sum_{m \leq x-1 \atop f(m)=-f(m+1)} 2B\ \leq\ 
{4B x \over L(x)^{7}},
$$
as claimed.  \qed
\end{pf} 

An immediate corollary of this Lemma is as follows:

\begin{cor} \label{shift_corollary}
If (\ref{main_consequence}) fails to hold, then for any integer $h \geq 1$ and
real number $\theta \leq 2$,
\begin{eqnarray}
\sum_{n \leq x/2 - h} |f(\lfloor n\theta\rfloor ) - 
f(h \lfloor n\theta/h \rfloor )| &\leq& {2 \over \theta} 
\sum_{m \leq x - 2h} |f(m) - f(h \lfloor m/h \rfloor ) | \nonumber \\
&\leq& {8h x \over \theta L(x)^{7}}. \nonumber
\end{eqnarray}
\end{cor}

\begin{pf}
We note that the first inequality of this lemma follows since for
each integer $m \geq 1$, there are at most $2/\theta$ integers $n$
such that $\lfloor n\theta \rfloor = m$ (and note that 
$\lfloor m/h \rfloor = \lfloor \lfloor n\theta \rfloor / h \rfloor = 
\lfloor n\theta / h \rfloor$); 
and, the second inequality is
an immediate consequence of Lemma \ref{shift_lemma}, since 
$h \lfloor n/h \rfloor = n + \delta(n)$ where $-h < \delta(n) \leq 0$.  \qed
\end{pf}

Finally, we have one more corollary:

\begin{cor} \label{shift_corollary2}
For $x$ sufficiently large, if (\ref{main_consequence}) fails to hold, 
and if $p,q \leq \beta_1$ are primes with $f(p)=f(q)$ and 
$p/q \in (1/2,2)$, then
$$
\Sigma(I,p/q)\ <\ {|I| \over \beta_2^2}.
$$
\end{cor}

\begin{pf}
To prove this corollary we have from the triangle inequality, from the above 
lemma and corollary, and from the fact that $f(p m) = f(p)f(m) = f(q)f(m) = f(qm)$
for $m \geq 1$,
\begin{eqnarray}
2\Sigma(I,p/q) &=& \sum_{n \in I} |f(n) - f(\lfloor pn/q \rfloor )| \nonumber \\
&\leq& \sum_{n \in I} |f(n) - f(q \lfloor n/q \rfloor)|\ +\ 
|f(q \lfloor n/q \rfloor) - f(p \lfloor n/q \rfloor )|\nonumber \\
&&\hskip0.5in +\ |f(\lfloor pn/q \rfloor) - f(p \lfloor n/q \rfloor ) |
\nonumber \\
&=& \sum_{n \in I} |f(n) - f(q\lfloor n/q \rfloor )| + 
|f(\lfloor pn/q \rfloor ) - f(p \lfloor n/q \rfloor)| \nonumber \\
&\leq& {4 q x \over L(x)^{7}}\ +\ {8q x \over L(x)^{7}}\ <\ {2|I| \over \beta_2^2},
\nonumber
\end{eqnarray}
for $x$ sufficiently large, which proves the corollary.  \qed 
\end{pf}

Suppose that $\alpha' \in R(x)$ is arbitrary; we may write it as
$$
\alpha' = {p_1 p_2 \cdots p_k \over q_1 q_2 \cdots q_k} = \alpha_1 
\alpha_2 \cdots \alpha_k,\ {\rm where\ }\alpha_i = {p_i \over q_i} \in 
\left ( 1/2,\ 2\right ),
$$
and where $k < \beta_2$ and $p_1,...,p_k,q_1,...,q_k < \beta_1$ 
are primes satisfying $f(p_i)=f(q_i)$.  
In addition, we may assume that
$\alpha_1\alpha_2\cdots \alpha_\ell \in (1/2,2)$ for all $1 \leq \ell \leq k$;
that is, through a simple induction argument, one can show that there exists 
some permutation $\sigma$ of $1,2,...,k$, such that
$\alpha_{\sigma(1)}\alpha_{\sigma(2)}\cdots \alpha_{\sigma(\ell)} 
\in (1/2,2)$, for $\ell=1,2,...,k$.  

We will show that our assumption that (\ref{main_consequence}) fails to
hold implies  
\begin{equation}\label{alpha_induction}
\Sigma(I, \alpha_1\alpha_2 \cdots \alpha_\ell)\ \leq\ {|I|\ell \over
\beta_2^2},\ {\rm for\ all\ } 1 \leq \ell \leq k.
\end{equation}
Such an inequality, if true, implies that (\ref{translation_inequality2})
fails to hold ( since $\alpha' = \alpha_1\cdots \alpha_k$, where
$k < \beta_2$); and so, since $\alpha' \in R(x)$ was arbitrary, we could
conclude that (\ref{translation_inequality2}) fails to hold for all 
$\alpha' \in R(x)$, which would prove (the contrapositive of) our Proposition. 

To see that (\ref{alpha_induction}) holds, we first note that it holds for 
$\ell=1$, since this follows from Corollary \ref{shift_corollary2}.
Now suppose that, for proof by induction, (\ref{alpha_induction}) 
holds for $\ell = m < k$.  Thus, 
$$
\sum_{n \in I} |f(n) - f(\lfloor \alpha_1\cdots \alpha_m n \rfloor )|
\ =\ 2\Sigma(I,\alpha_1\cdots \alpha_m)\ \leq\ {2|I| m \over \beta_2^2}. 
$$
We will show below that for $x$ sufficiently large, if $m \geq 1$, then
\begin{equation}\label{difference_inequality}
\sum_{n \in I} |f(\lfloor \alpha_1\cdots \alpha_m n \rfloor ) - f(\lfloor
\alpha_1\cdots \alpha_{m+1} n \rfloor ) | < {|I| \over \beta_2^2};
\end{equation}
and so, from this and our induction hypothesis we will have
\begin{eqnarray}
2\Sigma(I,\alpha_1\cdots \alpha_{m+1}) &=& \sum_{n \in I} 
|f(n) - f(\lfloor \alpha_1\cdots \alpha_{m+1} n\rfloor ) | \nonumber \\
&\leq& \sum_{n \in I} 
|f(\lfloor \alpha_1\cdots \alpha_m n \rfloor )-
f(\lfloor \alpha_1\cdots \alpha_{m+1} n \rfloor )| 
\nonumber \\
&&\hskip0.5in +\ 
|f(n) - f(\lfloor \alpha_1\cdots \alpha_m n \rfloor ) |
\nonumber \\ 
&\leq& {|I| \over \beta_2^2} + {2|I| m \over \beta_2^2}
\ <\ {2(m+1) |I| \over \beta_2^2}. \nonumber
\end{eqnarray}
which proves the induction step (that is, (\ref{alpha_induction}) holds
for $\ell=m+1$), and so (\ref{alpha_induction}) follows for 
$1 \leq \ell \leq k$.

Finally, we have from the triangle inequality and 
Corollary \ref{shift_corollary} that for $x$ sufficiently large, 
the left hand side of (\ref{difference_inequality}) is bounded
from above by
\begin{eqnarray}\label{f_string}
&& \sum_{n \in I} |f(\lfloor \alpha_1\cdots \alpha_m p_{m+1} n / q_{m+1} \rfloor )
- f(p_{m+1} \lfloor \alpha_1 \cdots \alpha_m n /q_{m+1}\rfloor )| \nonumber \\
&&\ \ \ \ \ +\ | f(\lfloor \alpha_1 \cdots \alpha_m n\rfloor ) - 
f( p_{m+1} \lfloor \alpha_1\cdots \alpha_m n/q_{m+1} \rfloor ) |\nonumber \\
&&\ \leq\ {8 q_{m+1} x \over \alpha_1\cdots \alpha_m L(x)^{7}}\nonumber \\
&&\hskip0.5in +\ \sum_{n \in I} | f(\lfloor \alpha_1 \cdots \alpha_m n \rfloor ) -
f( q_{m+1} \lfloor \alpha_1\cdots \alpha_m n/q_{m+1} \rfloor )|\nonumber \\
&&\ \leq\ {8 q_{m+1} x \over \alpha_1\cdots \alpha_m L(x)^{7}}
\ +\ {8 q_{m+1} x \over \alpha_1\cdots \alpha_m L(x)^{7}}\ \leq\ {|I| \over 
\beta_2^2}, \nonumber 
\end{eqnarray}
as claimed.   \qed
\bigskip

\section{Proof of Proposition \ref{prop_2}.}
\bigskip

Let $\beta_1, \beta_2, I, \Sigma(I,\alpha)$ and $R(x)$ be as in the 
Section \ref{main_theorem_section}.  We will need the following two results to prove
Proposition \ref{prop_2} (which are proved in Section \ref{lemma_section}): 

\begin{lem}\label{prop_2_smooth_lemma}
For $x$ sufficiently large and for any pair of real numbers $y_1,y_2$ satisfying
$$
1-{1 \over \beta_2}\ <\ y_1 < y_2 < 1,\ {\rm where\ }
y_2 - y_1 \geq {1 \over x^2},
$$
there exists $\theta \in R(x)$ with $\theta \in [y_1,y_2]$.
\end{lem}

\begin{lem}\label{prop_2_lemma_2}
For $x$ sufficiently large there exists a real number 
$\theta' \in (1 - \beta_2^{-1}, 1-(2\beta_2)^{-1})$, such that
$$
\Sigma(I,\theta')\ <\ {2|I| \over 3}.
$$
\end{lem}

To prove Proposition \ref{prop_2}, 
suppose that $x$ is sufficiently large so that we can apply Lemma \ref{prop_2_lemma_2},
and let $\theta'$ be as appears there.
If $\theta' = a/n$ for some integer $n \leq x$, then set $y_1 = \theta'$
and $y_2 = \theta'+1/x^2$; otherwise, if $\theta' \neq a/n$, for any integer
$n \leq x$ (and some integer $a$), then we let $y_1, y_2$ be such that
$\theta' \in [y_1,y_2]$, $y_2 - y_1 = 1/x^2$, and $[y_1,y_2]$ contains no
rationals of the form $a/n$, $n \leq x$.  We claim that for such $y_1$ and $y_2$,
we will have that for any $\kappa \in [y_1, y_2]$, 
\begin{equation}\label{equal_floor}
\lfloor y_1 n \rfloor = \lfloor \theta' n \rfloor = 
\lfloor \kappa n \rfloor = \lfloor y_2 n \rfloor,
\ {\rm for\ every\ }n \in I;
\end{equation}
and so, for any such $\kappa$, this gives 
\begin{equation} \label{equality_sigma}
\Sigma(I,y_1) = \Sigma(I,\theta') = \Sigma(I,\kappa) = \Sigma(I,y_2).
\end{equation}
Since $y_2 - y_1 = 1/x^2$, from Lemma \ref{prop_2_smooth_lemma} we have
that there exists $\theta \in R(x) \cap [y_1,y_2]$, and so
taking $\kappa = \theta$ in (\ref{equality_sigma}), we deduce
$$
\Sigma(I,\theta) =  \Sigma(I,\theta') < {2|I| \over 3}.\ \ \ \qed
$$

\section{Proofs of Lemmas.}\label{lemma_section}
\bigskip

\begin{pf*}{PROOF of Lemma \ref{prop_2_smooth_lemma}.} 

To prove the Lemma, we will construct a sequence of rational numbers
\begin{equation}\label{alpha_i_chain}
{1 \over 2}=\alpha_0 < \alpha_1 < \alpha_2 < \cdots 
\alpha_{t-1} < 1 - {1 \over 4x^2}< \alpha_t < \alpha_{t+1}=1,
\end{equation}
where 
\begin{equation}\label{alpha_i_inequality}
2 < {1 - \alpha_{i-1} \over 1 - \alpha_i}\ <\ {\beta_2 \over 8 \log^2 x}
\end{equation}
and where each $\alpha_i$ is of the form
$$
\alpha_i = {p_1 p_2 \cdots p_s \over q_1q_2\cdots q_s} \in R(x),
\ {\rm where\ }p_j{\rm 's}\ {\rm and\ }q_j{\rm 's\ are\ prime},
$$
and where $s < \log x$.  It is evident that $t < 4\log x$ for $x$ 
sufficiently large, since 
(\ref{alpha_i_chain}) and (\ref{alpha_i_inequality}) give us that
$$
4x^2 > {1 \over 1 - \alpha_{t-1}} > {2 \over 1 - \alpha_{t-2}} > \cdots 
> {2^{t-1} \over 1 - \alpha_0} = 2^t.
$$
 
If we had such a collection of $\alpha_i$'s, we claim that we could product 
together (the product can contain repeats) 
some of them, to produce a close rational approximation to $y_0 = (y_1+y_2)/2$;
moreover, this product will itself lie in $R(x)$.  To find such a product, we iterate
the following algorithm:
\bigskip

1.  Set $j=0$ and $n_0 = y_0$.

2.  Given the real number $n_j \in (1/2,1)$, 
we know that there exists $0 \leq i \leq t$ such that 
$\alpha_i \leq n_j < \alpha_{i+1}$.

3.  Set $n_{j+1} = n_j/\alpha_{i+1}$.

4.  Set $j \leftarrow j+1$, and repeat step 2 until 
$1 - 1/(4x^2) \leq n_j < 1$.
\bigskip

Let us assume for now that the algorithm terminates.  Then, if we let 
$\theta = \gamma_1\cdots \gamma_j$ be the product of all the $\alpha_i$'s 
we divide by each time the algorithm executes step 3, 
in passing from $n_0$ to $n_j$, we will have that
$$
1 - {1 \over 4x^2} \leq {n_0 \over \theta} = {y_0 \over \theta} < 1;
$$
and so, since $\theta = \gamma_1\cdots \gamma_j \leq 1$
we will have that $|y_0 - \theta| \leq 1/(4x^2)$, which gives us that
$\theta \in [y_1,y_2]$.  Now, in proving that
the algorithm above halts, we will show that $j < t \beta_2/(4\log^2 x)$;
and so, if we expand out $\theta$ in terms of a ratio of products
of primes (by expanding each $\gamma_i$ as such a product), we will have  
$$
\theta = {r_1\cdots r_g \over s_1\cdots s_g},\ r_i{\rm 's\ and\ } 
s_i{\rm 's\ are\ prime},
$$ 
and $r_i/s_i \in (1/2,2)$, $f(r_i)=f(s_i)$, and $g < j \log x < \beta_2$.
Thus, we will have $\theta \in R(x)$, and the Lemma will follow.

Let us now see that the above algorithm halts with 
$j < t \beta_2/(4\log^2 x) < \beta_2 / \log x$:  
We will show by induction that if $\alpha_i \leq n_0 < \alpha_{i+1}$ 
(where $i = 0,...,t$),
then the algorithm halts with $j < (t+1-i) \beta_2/(4\log^2 x)$.  
This clearly holds for $i=t$, since in this case, the above algorithm halts after the 
first pass, giving $\theta = 1 \in R(x)$.  
Suppose that, for proof by induction, the algorithm halts with 
$j < (t+1-i)\beta_2/(4\log^2 x)$ for all $n_0$ satisyfing 
$\alpha_i \leq n_0 < \alpha_{i+1}$, where $1 \leq B \leq i \leq t$.  

Now, If $n_0$ is such that $\alpha_{B-1} \leq n_0 < \alpha_B$, then the algorithm 
will repeatedly divide by $\alpha_B$ in step 3, say it divides $\ell$ times, 
until $n_\ell = n_0/\alpha_B^\ell \geq \alpha_B$.  Since the $\alpha_i$'s 
satisfy (\ref{alpha_i_inequality}), we will have that 
$\ell \leq \beta_2/(4\log^2 x)$.
We also have that $\alpha_B \leq n_\ell \leq 1$;
and therefore, $\alpha_{i'} \leq n_\ell < \alpha_{i'+1}$,
where $B \leq i' \leq t$.  By the induction hypothesis, we will have that
$j-\ell < (t+1-B)\beta_2/(4\log^2 x)$ 
(start the algorithm with $n_0$ equal to this number
$n_\ell$); so, the algorithm will halt with
$j < (t+1-(B-1))\beta_2/(4\log^2 x)$, 
since $\ell < \beta_2/(4\log^2 x)$, and so the induction step is proved.

To finish the proof of the Lemma, we must prove that the $\alpha_i$'s
exist:  Suppose that, for proof by induction, we have constructed the
numbers $\alpha_0,...,\alpha_u$.  Let $y = (1-\alpha_u)^{-1}(16\log^2 x)^{-1}
\beta_2$.  We will show that there exists integers $y < m_1 < m_2 < 2y$,
such that $1 \leq m_2 - m_1 < \beta_2(32 \log^2 x)^{-1}$, which will give
$$
2 < {1 - \alpha_u \over 1 - m_1/m_2} < {\beta_2 \over 8\log^2 x}.
$$
Moreover, $m_1$ and $m_2$ will each be expressible as 
$$
m_1 = p_1\cdots p_s,\ m_2 = q_1\cdots q_s,\ p_i\ {\rm and\ }q_i\ {\rm prime
\ for\ }i=1,...,s, 
$$
where $f(p_i)=f(q_i)$ and for $i=1,2,...s$,
$$
y^{1/s} < p_i, q_i < y^{1/s} \left ( 1 + {1 \over 2s} \right ),
{\rm where\ } s = \left \lfloor {\sqrt{\log y} \over 4} \right \rfloor. 
$$
One can easily check that for $x$ sufficiently large 
$s < \sqrt{\log x}$ (this follows from (\ref{alpha_i_chain}) and
(\ref{alpha_i_inequality}), which give that $(1-\alpha_u)^{-1} < x^3$ ), 
$y^{1/s} < \beta_1/2$, and
that $p_i/q_i \in (1/2,2)$.  Thus, $m_1/m_2 \in R(x)$ with
$s < \log x$ as claimed; and so, letting $\alpha_{u+1} = m_1/m_2$
will prove the induction step above, giving that the $\alpha_i$'s 
exist.

We are left to prove that $m_1$ and $m_2$ exist:  Let $S(y)$ denote the set
of all integers which can be written as $m = p_1\cdots p_s$,
where $p_i \in J := [y^{1/s}, y^{1/s}
(1 + (2s)^{-1})]$ is prime.
We note that all these primes are $\leq \beta_1$, and all 
such products lie in $[y,2y]$.  
Also, by the pigeonhole principle, there exists a subset 
$T(y) \subseteq S(y)$, with $|T(y)| \geq |S(y)|/s$,
and a constant $1 \leq D \leq s$, such that if 
$m \in T(y)$ has prime factorization $m = p_1\cdots p_s$,
then exactly $D$ of these prime factors $p_i$ (counting repeats) will
have $f(p_i) = 1$ (and so, $s-D$ factors $p_j$ will have $f(p_j)=-1$).
Thus, for any pair of numbers $m, m' \in T(y)$, we can arrange their
prime factorizations so that
$$
m = p_1\cdots p_s,\ m' = q_1\cdots q_s,\ f(p_i) = f(q_i).
$$ 
Thus, to prove that $m_1,m_2$ exist, it suffices to show that 
$|T(y)| > (32y\log^2 x)\beta_2^{-1}$ (since this implies there are
at least two elements of $T(y) \subseteq [y,2y]$ which are at most
$\beta_2 (32\log^2 x)^{-1}$ apart):  
By the Prime Number Theorem, for $x$ sufficiently large there
are $> y^{1/s}/(3\log y)$ primes in $J$; and so, by some elementary
combinatorics and the bounds $s < \sqrt{\log x}$,  $n! \leq n^{n-1}$, and
$\log y < 3 \log x$, we get
\begin{eqnarray}
|T(y)|\ &>&\ {|S(y)| \over s}\ >\ {(y^{1/s}/(3\log y) )^s
\over s!s}\ \geq\  y\left ( {1 \over 3s \log y} \right )^s 
\nonumber \\
&>&\ {y \over \exp(2\log\log x \sqrt{\log x})}\ >\ {32 y \log^2 x
\over \beta_2}. \nonumber \\
\end{eqnarray}
The Lemma now follows.   \qed
\end{pf*}

\begin{pf*}{PROOF of Lemma \ref{prop_2_lemma_2}.}

Let $u_0 = 1 - (2\beta_2)^{-1}$, $u_1 = 1 - \beta_2^{-1}$.
Since $\Sigma(I,t) \geq 0$ for all $t$ real, to prove the lemma it 
suffices to show that 
\begin{equation}\label{integral_inequality}
\int_{u_0}^{u_1} \Sigma(I, t) dt\ \sim\ {(u_1 - u_0)|I| \over 2}.
\end{equation}

We have
\begin{eqnarray}
\int_{u_0}^{u_1} \Sigma(I,t) dt &=& \int_{u_0}^{u_1} \#\{ n \in I\ :\ 
f(n) = -f(\lfloor t n \rfloor )\} dt \nonumber \\ 
&\leq& \sum_{n \in I} \sum_{u_0 n \leq m \leq u_1 n\atop f(m) = -f(n)} 
\mu\left ( \{ t\ :\ m = \lfloor t n \rfloor\} \right ), \nonumber 
\end{eqnarray}
where $\mu$ is the Lebesgue measure.

We have that $m = \lfloor t n \rfloor$ if and only if 
$m/n \leq t < (m+1)/n$; and so,
$$
\mu \left ( \{ t\ :\ m = \lfloor t n \rfloor \} \right ) = {1 \over n}.
$$
Thus,
\begin{equation}\label{reduce_measure}
\int_{u_0}^{u_1} \Sigma(I,t) dt \leq \sum_{n \in I} \sum_{u_0 n \leq m \leq u_1 n
\atop f(m) = -f(n)} {1 \over n}\nonumber
\end{equation}
To estimate this inner sum we use the following result of Hildebrand (see
Corollary 1 of \cite{hildebrand}):

\begin{thm}\label{hildebrand_theorem} 
Let $f(n)$ be a real-valued multiplicative function of modulus
$\leq 1$ satisfying (\ref{main_assumption}), and let $\phi(z)$ satisfy
$$
3 \leq \phi(z) \leq z,\ \log \phi(z) \sim \log z\ (z \to \infty).
$$
Then the limit
$$
\lim_{z \to \infty} {1 \over \phi(z)}\sum_{z-\phi(z) < n < z} f(n) = 0. 
$$
\end{thm}

From this Theorem, with $z=u_1n$ and $\phi(z) = n(u_1-u_0)$ we deduce
for $x\exp(-\log^{2/3}x) < n < x$ 
$$
\sum_{u_0 n \leq m \leq u_1 n\atop f(m)=-f(n)} 1\ \sim\ {n(u_1-u_0) \over 2};
$$
and so,
$$
\sum_{n \in I} {1 \over n} \sum_{u_0 n \leq m \leq u_1n\atop f(m)=-f(n)} 1 
\ \sim\ \sum_{n \in I} {1 \over n}\ {(u_1 - u_0)n \over 2}\ =\ 
{(u_1-u_0)|I| \over 2}.
$$
Putting this into (\ref{reduce_measure}), we deduce (\ref{integral_inequality}),
and our Lemma is proved.    \qed
\end{pf*}

\end{document}